\documentclass[11pt,reqno]{article}
\usepackage{hyperref}

\usepackage[light,first,bottomafter]{draftcopy}
\usepackage{amsmath}
\usepackage{amssymb}
\usepackage{amsthm}
\usepackage{srcltx}
\usepackage{enumerate}
\numberwithin{equation}{section}

\newtheorem{lemma}{Lemma}[section]
\newtheorem{theorem}[lemma]{Theorem}

\newtheorem{remarks}[lemma]{Remarks}

\def\2to{\to\!\!\!\!\!\!{}_{_2}\;\;}
 
\def\tauto{\to\!\!\!\!\!\!{}_{_\tau}\;\;}
 
\def\w2to{\rightharpoonup\!\!\!\!\!\!{}_{_2}\;\;}
\def\ws2to{\rightharpoonup\!\!\!\!\!\!{}_{_2}\!\!\!{}^{*}\;\;\;}
\def\wto{\rightharpoonup}

\def\erre
{{\bf R}}

\def\enne
{{\bf N}}

\def\Theta{{\mit \Theta}}
\let\La=\Lambda\def\Lambda{{\mit \La}}
\let\Sg=\Sigma\def\Sigma{{\mit \Sg}}
\let\Pig=\Pi\def\Pi{{\mit \Pig}}
\let\Ps=\Psi\def\Psi{{\mit \Ps}}
\let\Xii=\Xi\def\Xi{{\mit \Xii}}

\title{\bf Evolutionary {\boldmath $\Gamma$}-convergence of weak-type}  
\author{Augusto Visintin
\thanks{Dipartimento di Matematica dell'Universit\`a degli Studi di Trento --
via Sommarive 14,  38050 Povo di Trento, Italia -- email: augusto.visintin@unitn.it }
}

\date{\today} 
\begin{document}
\maketitle 
   
\begin{abstract}  
A notion of evolutionary $\Gamma$-convergence of weak type is introduced
for sequences of operators acting on time-dependent functions.
This extends the classical definition of $\Gamma$-convergence of functionals due to
De Giorgi. 
The $\Gamma$-compactness of equi-coercive and equi-bounded sequences of 
operators is proved.
Applications include the {\it structural compactness and stability\/} of quasilinear flows
for pseudo-monotone operators.  
\end{abstract}  

 
\bigskip 
\noindent{\bf AMS Classification (2000):}  
35K60, % Boundary value problems for nonlinear parabolic PDE  
47H05, % Monotone operators (with respect to duality)  
49J40, % Variational methods including variational inequalities      
58E. % Variational problems in infinite-dimensional spaces 

\section{Introduction}
\label{intro} 

\noindent 
In this note we deal with a notion of {\it evolutionary $\Gamma$-convergence of weak type,\/}
which extends De Giorgi's basic definition;
see \cite{DeFr} and e.g.\ the monographs \cite{At}, \cite{Bra1}, \cite{Bra2}, \cite{Da}. 

Here we formulate $\Gamma$-convergence for operators 
(rather than functionals) that act on time-dependent functions ranging in a Banach space $X$. 
This definition of evolutionary $\Gamma$-convergence is 
quite different from that of \cite{SaSe} 
and from that of \cite{DaSa}, \cite{Mi1}, \cite{Mi2}.
In those works $\Gamma$-convergence is set for almost any $t\in {}]0,T[$, 
whereas here it is assumed just in the weak topology of $L^1(0,T)$, 
and thus is a substantially weaker notion. 
This definition fits the rather general framework of 
$\bar\Gamma$-convergence, defined in Chap.~16 of \cite{Da}; see also references therein. 
However \cite{Da} does not encompass Theorem~\ref{teo.comp} ahead,
which is the main achievement of this work and is at the basis of the results of \cite{Vi17};
see also the survey \cite{ViLinc}. 

After results of Brezis and Ekeland \cite{BrEk}, Nayroles \cite{Na}
and Fitzpatrick \cite{Fi}, a class of first-order quasilinear flows can be given a variational formulation, without assuming monotonicity of the operator.
In the parallel work \cite{Vi17}, Theorem~\ref{teo.comp} is applied
to prove structural compactness and structural stability of the corresponding Cauchy problem. 
This applies to a large number of PDEs of mathematical physics, see \cite{Vi17}.

We proceed as follows. 
In Section~\ref{sec.evol} we define evolutionary $\Gamma$-convergence of weak type.
In Theorem~\ref{teo.comp} we prove compactness for this convergence.
In Theorem~\ref{teo.comp'} we apply Theorem~\ref{teo.comp} 
to the structural compactness and stability of nonmonotone flows.

\section{Definition of evolutionary {\boldmath $\Gamma$}-convergence of weak type}
\label{sec.evol}  

\noindent
In this section we extend De Giorgi's basic notion of $\Gamma$-convergence to operators 
(rather than functionals) that act on time-dependent functions ranging in a Banach space $X$.  
 
\bigskip
\noindent{\bf Functional set-up.}
Let $X$ be a real separable and reflexive Banach space, $p\in [1,+\infty[$, $T>0$, and 
\begin{equation}\label{eq.evol.mu}
\begin{split} 
&\hbox{ $\mu$ be a positive and finite measure on $]0,T[$, }
\\ 
&\hbox{ that is absolutely continuous w.r.t.\ the Lebesgue measure. } 
\end{split} 
\end{equation} 
Examples of interest will be the Lebesgue measure of $]0,T[$,
and $\mu$ such that $d\mu(t) = (T-t) \, dt$. 
Let us set
\begin{equation}\label{eq.evol.defgamma.1}
\begin{split}
&L^p_\mu(0,T) = \Big\{\mu\hbox{-measurable }v: {}]0,T[{}\to \erre:
\!\!\int_0^T \! |v(t)|^p \, d\mu(t)< +\infty \Big\},
\\
&L^p_\mu(0,T;X) = \Big\{\mu\hbox{-measurable }w: {}]0,T[{}\to X:
\!\!\int_0^T \! \|w(t)\|_X^p \, d\mu(t)< +\infty \Big\},
\end{split}
\end{equation} 
and equip either space with the respective graph norm.  

Let us equip $L^p_\mu(0,T;X)$ with a topology $\tau$ that is intermediate 
between the weak and the strong topology.
\footnote{ The generality of this topology is instrumental to the application 
to structural stability, see \cite{Vi17}. 
A reader interested just in evolutionary $\Gamma$-convergence might go through this section assuming that $\mu$ is the Lebesgue measure and that $\tau$ is the weak topology.
}
For any operator $\psi: L^p_\mu(0,T;X)\to L^1_\mu(0,T):w\mapsto \psi_w$, let us set 
\begin{equation}\label{eq.evol.crochet}
[\psi,\xi](w) = \int_0^T \psi_w(t) \, \xi(t) \, d\mu(t) 
\qquad\forall w\in L^p_\mu(0,T;X), \forall \xi\in L^\infty(0,T). 
\end{equation}  
 
\medskip
\noindent{\bf Evolutionary $\Gamma$-convergence of weak type.}
Let $\{\psi_n\}$ be a bounded sequence of operators $L^p_\mu(0,T;X)\to L^1_\mu(0,T)$; 
by this we mean that, for any bounded subset $A$ of $L^p_\mu(0,T;X)$, 
the set $\{\psi_{n,w}: w\in A, n\in \enne\}$ is bounded in $L^1_\mu(0,T)$.
If $\psi$ is also an operator $L^p_\mu(0,T;X)\to L^1_\mu(0,T)$, we shall say that 
\begin{equation}\label{eq.evol.defgamma.1+}
\begin{split}
&\hbox{ $\psi_n$ sequentially $\Gamma$-converges to }\psi
\\
&\hbox{ in the topology $\tau$ of $L^p_\mu(0,T;X)$ and } 
\\
&\hbox{ in the weak topology of }L^1_\mu(0,T)
\end{split}
\end{equation}
if and only if, denoting by $L^\infty_+(0,T)$ the cone of the nonnegative functions 
of $L^\infty(0,T)$,
\begin{equation}\label{eq.evol.defgamma.2}
[\psi_n, \xi] \mbox{ sequentially $\Gamma\tau$-converges to $[\psi,\xi]$ in }
L^p_\mu(0,T;X), \; \forall \xi\in L^\infty_+(0,T).
\end{equation}  
We shall say that a sequence {\it $\Gamma\tau$-converges\/} if it $\Gamma$-converges with respect to $\tau$, that a functional is {\it $\tau$-lower semicontinuous\/} if it is lower semicontinuous with respect to $\tau$, and so on.

By the classical definition of sequential $\Gamma$-convergence,  
\eqref{eq.evol.defgamma.2} means that for any $\xi\in L^\infty_+(0,T)$ 
\begin{equation}\label{eq.evol.defgamma.3}
\begin{split}
&\forall w\in L^p_\mu(0,T;X), \forall \hbox{ sequence $\{w_n\}$ in }L^p_\mu(0,T;X), 
\\
&\hbox{if \ $w_n\tauto w$ in $L^p_\mu(0,T;X)$ \ then \ }
\liminf_{n\to+ \infty} \; [\psi_n,\xi](w_n) \ge [\psi,\xi](w), 
\end{split}
\end{equation}
  \vskip-0.4truecm
\begin{equation}\label{eq.evol.defgamma.4}
\begin{split}
&\forall w\in L^p_\mu(0,T;X), \exists\hbox{ sequence $\{w_n\}$ 
of $L^p_\mu(0,T;X)$ such that } 
\\
&w_n\tauto w \hbox{ in $L^p_\mu(0,T;X)$ \ and \ }
\lim_{n\to +\infty} \; [\psi_n,\xi](w_n) = [\psi,\xi](w). 
\end{split}
\end{equation}

By the properties of ordinary $\Gamma$-convergence, \eqref{eq.evol.defgamma.2} entails that
\begin{equation}\label{eq.evol.defgamma.5}
[\psi, \xi] \hbox{ is sequentially $\tau$-lower semicontinuous in }L^p_\mu(0,T;X), \; 
\forall \xi\in L^\infty_+(0,T).
\end{equation}
 
\begin{remarks}\rm
(i) This definition of {\it evolutionary $\Gamma$-convergence\/} is not equivalent either 
to that of \cite{SaSe} or to that of \cite{DaSa}, \cite{Mi1}, \cite{Mi2}.
In those works $\Gamma$-convergence is actually assumed for almost any 
$t\in {}]0,T[$, whereas here it is just weak in $L^1_\mu(0,T)$. 

(ii) The present definition is based on testing the sequence $\{\psi_{n,w}\}$ on functions of time, 
and may equivalently be reformulated in terms of set-valued functions as follows.
Denoting by ${\cal L}(0,T)$ the $\sigma$-algebra of the 
Lebesgue-measurable subsets of $]0,T[$, 
\eqref{eq.evol.defgamma.2} is tantamount to
\begin{equation}\label{eq.evol.defgamma.2'}
\int_A \psi_{n,w}(t) \, d\mu(t) \to \int_A \psi_w(t) \, d\mu(t) 
\qquad\forall A\in {\cal L}(0,T).
\end{equation} 
As the elements of ${\cal L}(0,T)$ are in one-to-one correspondence 
with the characteristic functions of $L^\infty(0,T)$,
this equivalence can be checked by mimicking the argument based on the Lusin theorem,
that we use in the proof below.

(iii) By restating the definition \eqref{eq.evol.defgamma.1+}, \eqref{eq.evol.defgamma.2} as in the previous remark, it fits the rather general framework of $\bar\Gamma$-convergence, 
that is defined in Chap.~16 of \cite{Da}, see also references therein.
The results of that monograph however do not encompass the theorem of
$\Gamma$-compactness of the next section. 
  
(iv) Although we consider generic operators $\psi: L^p(0,T;X)\to L^1(0,T)$,
our main concern is for superposition (i.e., Nemytski\u{\i}-type) operators of the form  
\begin{equation}\label{eq.evol.super}
\begin{split}
&\psi_w(t) = \varphi(t,w(t))  
\qquad\forall w\in L^p_\mu(0,T;X),\hbox{ for a.e.\ }t\in {}]0,T[,
\\
&\varphi: {}]0,T[{} \times X\to \erre^+ \hbox{ \ being a {\it normal function.\/} }
\end{split}
\end{equation} 
(By this we mean that $\varphi$ is globally measurable and $\varphi(t,\cdot)$ is lower semicontinuous for a.e.\ $t\in {}]0,T[$.)
\hfill$\Box$  
\end{remarks}

\section{Evolutionary {\boldmath $\Gamma$}-compactness of weak type}
\label{sec.comp}  

\noindent
In this section we prove a theorem of compactness for the evolutionary $\Gamma$-convergence
that we just defined.  
This will be based on a result of Hiai \cite{Hi}, that we now display.

As a preparation, let us say that a functional 
$\Phi: L^p_\mu(0,T;X)\to \erre$ is invariant by translations if,
denoting by $\widetilde w$ the function obtained by extending $w$ to $\erre$ 
with null value, 
\begin{equation}\label{eq.evol.invar}
\begin{split} 
&\Phi\big(\widetilde w(\cdot +s)\big|_{[0,T]}\big) = \Phi(w)
\\[1mm]
&\forall w\in L^p_\mu(0,T;X), \forall s>0 
\text{ such that $w=0$ a.e.\ in }]s,T[. 
\end{split}
\end{equation} 

The next result will play a role in the present analysis.

\begin{lemma} [\cite{Hi}] \label{lemma.Hiai}
Let $X$ be a real separable Banach space and $p\in [1,+\infty[$.
Let a functional $\Phi: L^p_\mu(0,T;X)\to \erre\cup \{+\infty\}$ 
($\Phi \not\equiv +\infty$) be lower semicontinuous and also {\rm additive,\/} 
in the sense that  
\begin{equation}\label{eq.evol.add}
\begin{split} 
&\forall w_1,w_2\in L^p_\mu(0,T;X), 
\\
&\mu\big( \{t\in {}]0,T[{}: w_1(t)w_2(t) \not= 0\} \big) =0
\quad\Rightarrow\quad \Phi(w_1 +w_2) = \Phi(w_1) + \Phi(w_2).
\end{split}
\end{equation}
\indent
Then there exists a normal function $\varphi: {}]0,T[{} \times X\to \erre\cup\{+\infty\}$ 
such that
\begin{eqnarray}
\begin{split}\label{aaa}
&\varphi(t,\cdot) \not\equiv +\infty \quad\forall t\in ]0,T[,
\\
&\varphi(\cdot,0) =0 \quad\hbox{ a.e.\ in }{}]0,T[,
\\ 
&\Phi(w) = \int_0^T \varphi(t, w(t)) \, d\mu(t)
\qquad\forall w\in L^p_\mu(0,T;X).
\end{split} 
\end{eqnarray}  
\indent
Moreover, if $\Phi$ is convex then $\varphi(t,\cdot)$ is also convex for a.e.\ $t\in{}]0,T[$,
and the function $\varphi$ is unique, up to modification on sets of the form 
$N\times X$ with $\mu(N)=0$.
\\
\indent
Finally, if the functional $\Phi$ is invariant by translations, 
then $\varphi$ does not depend on $t$.
\end{lemma}  

This lemma may be compared e.g.\ with Section~2.4 of \cite{But},
which however deals with a finite-dimensional space $X$. 

We are now ready to state and prove the main result of this note.
  
\begin{theorem} \label{teo.comp}
Let $X$ be a real separable and reflexive Banach space, $p\in [1,+\infty[$, $T>0$,
and $\{\varphi_n\}$ be a sequence of normal functions 
${}]0,T[{} \times X\to \erre^+$.
Assume that this sequence is equi-coercive and equi-bounded, in the sense that  
\begin{equation}\label{eq.evol.equibc} 
\begin{split}
&\exists C_1,C_2,C_3 >0: \forall n,\hbox{for a.e.\ }t\in {}]0,T[, \forall w\in X, 
\\
&C_1 \|w\|_X^p\le \varphi_n(t,w) \le C_2\|w\|_X^p +C_3, 
\end{split} 
\end{equation}
and that
\begin{equation}\label{eq.evol.nul}
\varphi_n(t,0) =0 \qquad\hbox{ for a.e.\ }t\in {}]0,T[,\forall n.
\end{equation}  
\indent
Let $\mu$ fulfill \eqref{eq.evol.mu}. Let
$\tau$ be a topology on $L^p_\mu(0,T;X)$ that either coincides or is finer 
than the weak topology, and such that 
\begin{equation}\label{eq.evol.gamcom} 
\begin{split}
&\hbox{ for any sequence $\{F_n\}$ of functionals }
L^p_\mu(0,T;X)\to \erre^+\cup \{+\infty\},
\\
&\hbox{ if \ }\sup_{n\in\enne} \big\{\|w\|_{L^p_\mu(0,T;X)}: 
w\in L^p_\mu(0,T;X), F_n(w) \le C \big\}<+\infty,
\\
&\hbox{ then \ $\{F_n\}$ has a sequentially $\Gamma\tau$-convergent subsequence. }  
\end{split}
\end{equation}
\indent
Then there exists a normal function $\varphi: {}]0,T[{} \times X\to \erre^+$ such that
$\varphi(\cdot,0) =0$ a.e.\ in $]0,T[$, and such that, 
defining the operators $\psi,\psi_n: L^2_\mu(0,T;X)\to L^1_\mu(0,T)$ for any $n$ 
as in \eqref{eq.evol.super}, possibly extracting a subsequence,
\begin{equation}\label{eq.evol.tesi} 
\begin{split}
&\hbox{ $\psi_n$ sequentially $\Gamma$-converges to }\psi
\\
&\hbox{ in the topology $\tau$ of $L^p_\mu(0,T;X)$ and } 
\\
&\hbox{ in the weak topology of $L^1_\mu(0,T)$ (cf.\ \eqref{eq.evol.defgamma.1+}). }
\end{split}
\end{equation}
\indent
Moreover, if $\varphi_n$ does not depend on $t$ for any $n$, 
then the same holds for $\varphi$.
\end{theorem}   

\noindent{\bf Proof.\/}  
For the reader's convenience, we split this argument into several steps.

(i) First we show that, denoting by $C^0_+([0,T])$ the cone of the nonnegative functions of 
$C^0([0,T])$, 
\begin{equation}\label{eq.evol.conv}
\begin{split}
&\forall \xi\in C^0_+([0,T]), \exists g_\xi: L^p_\mu(0,T;X)\to \erre^+ 
\hbox{ such that } 
\\
&[\psi_n,\xi] \hbox{ $\Gamma\tau$-converges to }g_\xi.
\end{split}
\end{equation}

The separable Banach space $C^0([0,T])$ has a countable dense subset $M$,
e.g., the family of polynomials with rational coefficients. Let us
denote by $M_+$ the cone of the nonnegative elements of $M$.
For any $\xi\in M_+$, by \eqref{eq.evol.equibc} a suitable subsequence $\{[\psi_{n'},\xi]\}$ 
weakly $\Gamma$-converges to a function $g_\xi: L^p_\mu(0,T;X)\to \erre$, and
\begin{equation}\label{eq.evol.a}
\begin{split}
&C_1 \int_0^T \|w(t)\|_X^p \, \xi(t) \, d\mu(t) \le g_\xi(w) 
\\
&\le C_2\int_0^T \|w(t)\|_X^p \, \xi(t) \, d\mu(t) + C_3 \int_0^T \xi(t) \, d\mu(t)
\qquad\forall w\in L^p_\mu(0,T;X). 
\end{split}
\end{equation}  

A priori the selected subsequence $\{[\psi_{n'},\xi]\}$ might depend on $\xi$. 
However, because of the countability of $M_+$, via a diagonalization procedure 
one can select a subsequence that is independent of $\xi\in M_+$.
(Henceforth we shall write $\psi_n$ in place of $\psi_{n'}$, dropping the prime.) 
For that subsequence thus
\begin{equation}\label{eq.evol.a=}
\begin{split}
&\forall \xi\in M_+, \exists g_\xi: L^p_\mu(0,T;X)\to \erre^+ \hbox{ such that } 
\\
&[\psi_n,\xi]\hbox{ \ $\Gamma\tau$-converges to \ $g_\xi$ in }L^p_\mu(0,T;X),
\end{split}
\end{equation}
that is, for any $\xi\in M_+$,
\begin{equation}
\begin{split}\label{eq.evol.b=}
&\forall w\in L^p_\mu(0,T;X), \forall \hbox{ sequence $\{w_n\}$ in }L^p_\mu(0,T;X), 
\\
&\hbox{if \ $w_n \tauto w$ in $L^p_\mu(0,T;X)$ \ then \ }
\liminf_{n\to+ \infty} \; [\psi_n,\xi](w_n) \ge g_\xi(w), 
\end{split}
\end{equation}
       \vskip-0.4truecm 
\begin{equation}
\begin{split}\label{eq.evol.c=}
&\forall w\in L^p_\mu(0,T;X), \exists\hbox{ sequence $\{w_n\}$ of $L^p_\mu(0,T;X)$ 
such that } 
\\
&w_n \tauto w \hbox{ in $L^p_\mu(0,T;X)$ \ and \ }
\lim_{n\to + \infty} \; [\psi_n,\xi](w_n) = g_\xi(w). 
\end{split}
\end{equation} 
(The recovery sequence $\{w_n\}$ in \eqref{eq.evol.c=} may depend on $\xi$.)

As any $\xi\in C^0_+([0,T])$ is the uniform limit of some sequence $\{\xi_m\}$ in $M$,
for any bounded sequence $\{w_n\}$ in $L^p_\mu(0,T;X)$
\begin{equation}
\begin{split}
&\sup_n \big| [\psi_n,\xi](w_n) - [\psi_n,\xi_m](w_n) \big| 
\overset{\eqref{eq.evol.crochet}}{=}
\sup_n \Big|\int_0^T \psi_{n,w_n}(t) \, [\xi(t) - \xi_m(t)] \, d\mu(t) \Big| 
\\
&\overset{\eqref{eq.evol.super}}{\le}
\| \xi - \xi_m \|_{C^0([0,T])} \sup_n \int_0^T \varphi_n(t,w) \, d\mu(t) 
\\
&\overset{\eqref{eq.evol.equibc}}{\le}\| \xi - \xi_m \|_{C^0([0,T])}
\sup_n \bigg\{ C_2 \!\! \int_0^T \|w_n(t)\|_X^p \, d\mu(t) 
+ C_3 \mu(]0,T[) \bigg\}
\qquad\forall m.
\end{split}
\end{equation}
By the density of $M_+$ in $C^0_+([0,T])$,  \eqref{eq.evol.conv} then follows.

\medskip
(ii) Next we extend \eqref{eq.evol.conv} to any $\xi\in L^\infty_+(0,T)$.

By the classical Lusin theorem, for any $\xi\in L^\infty_+(0,T)$ there exists a sequence 
$\{\xi_m\}$ in $C^0_+([0,T])$ such that,
setting $A_m =\{t\in [0,T]: \xi_m(t) \not=  \xi(t)\}$,
\begin{equation}\label{eq.evol.Lusin} 
\begin{split}
&\|\xi_m\|_{C^0([0,T])} \le \|\xi\|_{L^\infty(0,T)} \quad\forall m, 
\\
&\mu(A_m)\to 0.
\end{split}
\end{equation}  
Hence 
\begin{equation}
\begin{split}
&\sup_n \big| [\psi_n,\xi](w_n) - [\psi_n,\xi_m](w_n) \big| 
\overset{\eqref{eq.evol.crochet}}{=} 
\sup_n \Big|\int_0^T \psi_{n,w_n}(t) \, [\xi(t) - \xi_m(t)] \, d\mu(t) \Big| 
\\ 
&\overset{\eqref{eq.evol.super},\eqref{eq.evol.equibc}}{\le} \|\xi - \xi_m\|_{L^\infty(0,T)} 
\bigg\{ C_2 \sup_n \!\! \int_{A_m} \|w_n(t)\|_X^p \, d\mu(t) + C_3 \mu(A_m) \bigg\}
\qquad\forall m.
\end{split}
\end{equation}
As $w_n\tauto w$ in $L^p_\mu(0,T;X)$ (see \eqref{eq.evol.c=}), 
the sequence $\{\|w_n(\cdot)\|_X^p\}$ is equi-integrable. 
By this property and by \eqref{eq.evol.Lusin}$_2$, 
\[
\sup_n \int_{A_m} \|w_n(t)\|_X^p \, d\mu(t)\to 0
\qquad\text{ as }m\to \infty.
\]
\eqref{eq.evol.conv} is thus extended to any $\xi\in L^\infty_+(0,T)$. 

\medskip
(iii) Next we prove that  
\begin{equation}\label{eq.evol.due}
\begin{split}
&\exists \psi: L^p_\mu(0,T;X)\to L^1_\mu(0,T) \hbox{ such that } 
\\
&g_\xi(w) = [\psi_w,\xi] \qquad\forall w\in L^p_\mu(0,T;X),\forall \xi\in L^\infty_+(0,T).
\end{split}
\end{equation} 
(Here some care is needed, since $L^\infty_+(0,T)$ is not a linear space.)
Let us fix any $\xi\in L^\infty_+([0,T])$, any $w\in L^p_\mu(0,T;X)$, and any sequence
$\{w_n\}$ as in \eqref{eq.evol.c=}. 
By the boundedness of $\{w_n\}$ and by \eqref{eq.evol.equibc}, 
the sequence $\{\psi_{n,w_n}\}$ = $\{\varphi_n(\cdot,w_n)\}$ 
is bounded in $L^1_\mu(0,T)$ and is equi-integrable. 
There exists then a function $\gamma \in L^1_\mu(0,T)$ such that, 
possibly extracting a subsequence,
$\psi_{n,w_n}\wto \gamma$ in $L^1_\mu(0,T)$. 
\footnote{ We denote the strong and weak convergence respectively by $\to$ and $\wto$. 
} 
Thus
\begin{equation}\label{eq.evol.o}
[\psi_n,\xi](w_n) =\int_0^T \psi_{n,w_n}(t) \, \xi(t) \, d\mu(t) 
\to \int_0^T \gamma(t) \, \xi(t) \, d\mu(t)
\qquad\forall \xi\in L^\infty(0,T).
\end{equation}  
Thus by \eqref{eq.evol.c=}
\begin{equation}\label{eq.evol.q}
g_\xi(w) = \int_0^T \gamma(t) \, \xi(t) \, d\mu(t)
\qquad\forall \xi\in L^\infty_+(0,T).
\end{equation}

Therefore $\gamma$ is determined by $w\in L^p_\mu(0,T;X)$,  
and is independent of the specific sequence $\{w_n\}$ that fulfills \eqref{eq.evol.c=}.  
This defines an operator 
\begin{equation}\label{eq.evol.f}
\psi: L^p_\mu(0,T;X)\to L^1_\mu(0,T): w \mapsto \psi_w = \gamma.
\end{equation} 
The equality \eqref{eq.evol.q} thus reads
\begin{equation}\label{eq.evol.r}
g_\xi(w) =\int_0^T \psi_w(t) \, \xi(t) \, d\mu(t) 
\qquad\forall \xi\in L^\infty_+(0,T), \forall w\in L^p_\mu(0,T;X).
\end{equation} 
Recalling the definition \eqref{eq.evol.crochet}, 
we see that this completes the proof of \eqref{eq.evol.due}.  

\medskip
(iv) Finally, we show that there exists a normal function 
$\varphi: {}]0,T[{} \times X\to \erre^+$ such that the operator
$\psi$ that we just defined in \eqref{eq.evol.f} is as in \eqref{eq.evol.super}. 

By \eqref{eq.evol.nul}, for any $n$ the functional
\begin{equation}
\begin{split}
&\Phi_n: L^p_\mu(0,T;X)\to \erre^+: 
\\
&w\mapsto \int_0^T \psi_{w,n}(t) \, d\mu(t) =\int_0^T \varphi_n(t,w_n(t)) \, d\mu(t)
\end{split}
\end{equation}
is additive in the sense of \eqref{eq.evol.add} below. 
This property then also holds for the limit functional
\begin{equation}
\Phi: L^p_\mu(0,T;X)\to \erre^+: w\mapsto \int_0^T \psi_w(t) \, d\mu(t).
\end{equation}

By selecting $\xi\equiv 1$ in \eqref{eq.evol.defgamma.5}, we get that $\Phi$ is 
lower semicontinuous. 
By Lemma~\ref{lemma.Hiai} then there exists a normal function $\varphi$ as we just 
specified.
\hfill$\Box$

\section{Applications}
\label{sec.appl}  

\noindent
In this section we briefly illustrate how the notion of evolutionary $\Gamma$-convergence 
of weak type can be applied to prove the structural compactness and structural stability of flows 
of the form 
\begin{equation}\label{eq.flow1} 
D_tu + \alpha(u)\ni h 
\qquad\hbox{ in $V'$, a.e.\ in time }(D_t := {\partial/\partial t});
\end{equation} 
here $V$ is a Hilbert space, and
$\alpha:V\to {\cal P}(V')$ is a {\it semi-monotone\/} operator.
This is a particular case of the class of generalized pseudo-monotone operators of 
Browder and Hess \cite{BrHe}, and includes mappings of the form
\begin{equation}\label{eq.equicoer}
H^1_0(\Omega)\to {\cal P}(H^{-1}(\Omega)): 
v\mapsto - \nabla \cdot \vec\gamma(v,\nabla v),
\end{equation}
with $\vec\gamma$ continuous w.r..t. the first argument, and maximal monotone w.r..t. the
second one.
We refer to \cite{ViLinc} for a more expanded outline and to \cite{Vi17} 
for a detailed presentation. 

Under suitable restrictions, there exists a topology $\tau$ as above 
in Section~\ref{sec.evol}, such that
\begin{equation}\label{eq.fitzp.nonconv}
\begin{split} 
&\varphi:V \!\times\! V' \to \erre\cup \{+\infty\} 
\text{ is lower semicontinuous w.r.t.\ }\tau,
\\
&\varphi(v,v^*) \ge \langle v^*,v\rangle
\qquad\forall (v,v^*)\in V \!\times\! V',
\\
&\varphi(v,v^*) = \langle v^*,v\rangle
\quad\Leftrightarrow\quad v^*\in \alpha(v).
\end{split}
\end{equation} 
After defining the functional 
\begin{equation}
\begin{split} \label{eq.BENfun}
\Phi(v,v^*) = \int_0^T [\varphi(v,v^* -D_tv) - \langle v^*,v\rangle] \, d\mu(t) 
+ {1\over2}\int_0^T \|v(T)\|_H^2 \, dt - {T\over2} \|u(0)\|_H^2 &
\\
\forall (v,v^*)\in L^2(0,T;V \!\times\! V'),&
\end{split}
\end{equation}
one can show that
\begin{equation}\label{eq.equiv}
\eqref{eq.flow1} 
\quad\Leftrightarrow\quad 
\Phi(u,u^*) = \inf \Phi =0.
\end{equation}

This provides a (nonstandard) variational structure of the flow,
and paves the way to the use of a notion of evolutionary $\Gamma$-convergence
of weak type.
 
\bigskip 
\noindent{\bf Structural compactness and structural stability.}
We define {\it structural stability\/} 
as robustness to perturbations of the structure of the problem,
e.g.\ operators in differential equations.
These notions have obvious applicative motivations, as
data and operators are accessible just with some approximation.  
See e.g.\ \cite{ViLinc} and \cite{Vi17}.

Let us briefly illustrate these notions for a problem of the form $Au\ni h$, 
$A$ being a multi-valued operator acting in a Banach space and $h$ a datum.
Given bounded families $\{h_n\}$ and $\{A_n\}$, we formulate the stability of the problem
via two properties:

(i) {\it structural compactness:\/} existence of convergent sequences 
of data $\{h_n\}$ and of operators $\{A_n\}$ (in a sense to be specified);

(ii) {\it structural stability:\/} if $A_nu_n\ni h_n$ for any $n$, 
$A_n\to A$, $h_n\to h$ and $u_n\to u$, then $u$ is a solution of the asymptotic problem: 
$Au\ni h$.  

Structural compactness and structural stability of minimization principles 
can adequately be dealt with via De Giorgi's theory of $\Gamma$-convergence.
Next we outline how this can be extended to flows,
after these have been variationally formulated as in \eqref{eq.BENfun} and \eqref{eq.equiv}. 
  
\begin{theorem} [\cite{Vi17}] \label{teo.comp'}
Let $V$ be a real separable Hilbert space,
and $\mu$ be the measure on $]0,T[$ that fulfills \eqref{eq.evol.mu}.
Let $\{\varphi_n\}$ be a sequence of normal functions 
${}]0,T[{} \times V \!\times\! V'\to \erre^+$ such that  
\begin{eqnarray} 
&\varphi_n(t,\cdot) \in {\cal F}(V)
\qquad\hbox{for a.e.\ }t\in {}]0,T[, \forall n,  
\label{eq.comstab.reppn}
\\
&\begin{split}
&\exists C_1,C_2,C_3 >0: \forall n,\hbox{for a.e.\ }t\in {}]0,T[, \forall w\in V \!\times\! V', 
\\
&C_1 \|w\|_{V \!\times\! V'}^2\le \varphi_n(t,w) \le C_2\|w\|_{V \!\times\! V'}^2 +C_3, 
\end{split}   
\label{eq.comstab.equibc'}  
\\
&\varphi_n(t,0) =0 \qquad\hbox{ for a.e.\ }t\in {}]0,T[,\forall n,
\label{eq.comstab.nul'}
\end{eqnarray} 
and define the operators $\psi_n: L^2_\mu(0,T;V \!\times\! V')\to L^1_\mu(0,T)$ by
\begin{equation}\label{eq.comstab.super'} 
\psi_{n,w}(t) = \varphi_n(t,w(t))  
\qquad\forall w\in L^2_\mu(0,T;V \!\times\! V'),\hbox{ for a.e.\ }t\in {}]0,T[,\forall n.
\end{equation}  
\indent
Then there exists a normal function $\varphi: {}]0,T[{} \times V \!\times\! V'\to \erre^+$ 
such that
\begin{equation}\label{eq.comstab.repp}
\varphi(t,\cdot) \in {\cal F}(V)
\qquad\hbox{for a.e.\ }t\in {}]0,T[,
\end{equation}  
and such that, defining the corresponding operator 
$\psi: L^2_\mu(0,T;V \!\times\! V')\to L^1_\mu(0,T)$ 
as in \eqref{eq.comstab.super'}, possibly extracting a subsequence
\begin{equation}\label{eq.comstab.tesi'} 
\begin{split}
&\hbox{ $\psi_n$ sequentially $\Gamma$-converges to }\psi
\\
&\hbox{ in the topology $\widetilde\pi$ of $L^2_\mu(0,T;V \!\times\! V')$ and } 
\\
&\hbox{ in the weak topology of $L^1_\mu(0,T)$ (cf.\ \eqref{eq.evol.defgamma.1+}). }
\end{split}
\end{equation}
\indent
Moreover, if $\varphi_n$ does not depend on $t$ for any $n$, 
then the same holds for $\varphi$.
\end{theorem} 

This result provides the structural compactness and structural stability of flows of the form 
\begin{equation}\label{eq.flow2} 
D_tu - \nabla \cdot \vec\gamma(v,\nabla v) \ni h 
\qquad\hbox{ in $H^{-1}(\Omega)$, a.e.\ in time,}
\end{equation}  
with $\vec\gamma$ as above, see \cite{Vi17}. 
This can also be extended to doubly-nonlinear flows, see \cite{ViLinc}.

\bigskip
\centerline{\bf Acknowledgment}
\medskip

The author is a member of GNAMPA of INdAM.

This research was partially supported by a MIUR-PRIN 2015 grant for the project 
``Calcolo delle Variazioni" (Protocollo 2015PA5MP7-004).

This author is indebted to Giuseppe Buttazzo, who brought Hiai's paper \cite{Hi} to his attention.

\baselineskip=12truept
\bigskip

Author's address: 
\medskip

Augusto Visintin \par 
Universit\`a degli Studi di Trento \par 
Dipartimento di Matematica \par 
via Sommarive 14, \ 38050 Povo (Trento) - Italia \par 
Tel   +39-0461-281635 (office), +39-0461-281508 (secretary) \par 
Fax      +39-0461-281624 \par 
Email:   augusto.visintin@unitn.it 

\end{document}